\documentclass[twocolumn]{cinc}
\usepackage{graphicx}
\begin{document}
\bibliographystyle{cinc}

\title{Comparison of Two Formulations for Computing Body
Surface Potential Maps}


\author {Emma Lagracie$^{1}$, Lisl Weynans$^{1}$, Yves Coudière$^{1}$\\
\ \\ 
 $^1$ Univ. Bordeaux, CNRS, INRIA, Bordeaux INP, IMB, UMR 5251, IHU Liryc, F-33400 Talence, France}

\maketitle

\begin{abstract}
    In order to establish a link between the activation map and the body surface signals, we want to compute body potentials from activation maps and a frontlike approximation of the activation of the heart. To do so, two formulations that map the extracardiac potential from the transmembrane voltage naturally emerge from the bidomain model. Either the extracellular/extracardiac potential solves a Laplace equation with discontinuous conductivity coefficient and ionic current as a source (Source Formulation F1); or the quasi-stationary electrical balance between the intra- and extracellular fields (Balance Formulation F2). In this work, we compare F1 and F2, to determine which formulation is the most relevant to use.
    
    We compute reference activation map $\psi$, transmembrane voltage $v$ and body surface map $u$ with a bidomain 2D code. We design two alternative shapes $\tilde{v}$ for a frontlike approximation of $v$. Afterwards, two extracellular / extracardiac potentials are computed from the activation time $\psi$ and $\tilde{v}$, using the two different formulations. Then the extracardiac potentials solutions of F1 and F2 are compared respectively to the solution of the bidomain $u$.
    
    Results show that the Balance Formulation F2 is robust to the input data ($\tilde{v}$ and $\psi$). On the contrary, the Source Formulation F1 is very unstable and generates very large errors on the body surface map.

\end{abstract}

\section{Introduction}

In the framework of ECGi, we would like to reconstruct activation maps of the heart 
directly from torso potentials, available thanks to multi-electrode vests for example. 
To this end, we need to express a direct problem given the activation time $\psi(x)$. In particular, we want to deduce the extracellular and extracardiac potential $u$ directly from $\psi$.
We choose to model the transmembrane voltage as an activation front with
a predefined shape as in \cite{franzone1990mathematical, van2009non, ravon2019impact}. In this context, the approximate transmembrane voltage writes $\tilde{v}(x, t) = \bar{V}(t-\psi(x))$.
Even if its shape is uncertain, this representation of the transmembrane voltage captures the global behaviour of the depolarisation front. 

Then, the extracellular and extracardiac potential is expressed as a solution of a partial differential equation with a source term that depends on $\Bar{V}$, $\psi$, and the formulation of the operator $\left [\tilde{v} \longrightarrow u   \right ]$. This operator is derived from the bidomain equations that follow. In a heart ($\Omega_H$) embedded in a torso ($\Omega_T$) geometry, the bidomain model is given by
\begin{equation}
    \begin{aligned}
        &\text{div}(\sigma_i \nabla u_i) = \partial_t v  + I_{ion}(v) & \text{in } \Omega_H,\\
        &\text{div}(\sigma_e \nabla u) = -\partial_t v  - I_{ion}(v) & \text{in } \Omega_H,\\
        &\text{div}(\sigma_T \nabla u) = 0 & \text{in } \Omega_T,\\
        &\sigma_i \nabla (u+v) \cdot n = 0 & \text{on } \partial \Omega_H,\\
        & \sigma_T \nabla u \cdot n  = 0 & \text{on } \partial \Omega_T,
    \end{aligned}
    \label{bidomainF1}
\end{equation} with $u_i \in H^1(\Omega_H) ,~ u \in H^1(\Omega_H \cup \Omega_T)$ the intracellular and extracellular / extracardiac potentials respectively, $ v = u_i - u \in \Omega_H$ the transmembrane voltage, and $\sigma_i,~ \sigma_e,~ \sigma_T $ the conductivities respectively in the cardiac intracellular medium, cardiac extracellular medium and torso medium. They are scaled to use nondimensionnal terms. The ionic current, $I_{ion}$ is a non-linear (often cubic) function of $v$ and other variables that depend on the ionic model chosen. Equivalently, the volumic equations in the heart can write 
\begin{equation}
    \begin{aligned}
        &\text{div}(\sigma_i \nabla (u + v)) = \partial_t v  + I_{ion}(v) & \text{in } \Omega_H,\\
        &\text{div}((\sigma_i+\sigma_e) \nabla u) = -\text{div}(\sigma_i \nabla v) & \text{in } \Omega_H.
    \end{aligned}
    \label{bidomainF2}
\end{equation}
This version is more widely used in the literature. In these two equivalent formulations of the bidomain we see two equations modelling the dependency between $u$ and $v$: 
\begin{equation}
\label{SF}
\tag{F1}
    \begin{aligned}
  &-\text{div} \left(\sigma_e \nabla u \right) = \partial_t v + f(v) && \text{in } \Omega_H,\\
  &-\text{div} \left(\sigma_T \nabla u \right) = 0 && \text{in } \Omega_T;
\end{aligned}
\end{equation} and 
\begin{equation}
\label{BF}
\tag{F2}
    \begin{aligned}
      &- \text{div} \left(\sigma_e \nabla u \right) - \text{div} \left(\sigma_i \nabla (u+v) \right) =0 && \text{in } \Omega_H,\\
      &-\text{div} \left(\sigma_T \nabla u \right) = 0 && \text{in } \Omega_T.
    \end{aligned}
\end{equation}
Thus, \eqref{SF} and \eqref{BF} can be seen as two independent problems which data is $v$ and that define two mappings of $\left [ \tilde{v} \longrightarrow u \right ] $, taking $v=\tilde{v}$.
Either $u$ solves a Laplace equation with discontinuous conductivity coefficients (heart and torso) and a nonlinear source (Source Formulation \eqref{SF}), or it solves the quasi-stationary electrical balance between the intra and extracellular fields (Balance Formulation \eqref{BF}). Commonly, the potential $u$ is computed from the Balance Formulation \eqref{BF}. Anyway, we may also use \eqref{SF}, that has the advantage of not including any coupling on the heart boundary, and which differential operator is naturally shaped for the coupled heart/torso problem. Note that if $v$ solves the complete bidomain equations, the two formulations are equivalent. However, as soon as $v$ differs from the bidomain solution, \eqref{SF} and \eqref{BF} are no longer equivalent. In this study, we compare the solutions of the two formulations ($u_1$ for \eqref{SF} and $u_2$ for \eqref{BF}) evaluated in $\tilde{v}(x, t) = \bar{V}(t-\psi(x))$, with the bidomain solution (denoted by the index $_\text{ref}$). We also investigate the robustness of those two formulations to small variations in the transmembrane voltage and the activation time.

\section{Method}
\subsection{Reference solution and activation map} 
Let ($u_{\text{ref}}$, $v_{\text{ref}}$) be the couple solution of the bidomain with the Mitchell-Shaeffer \cite{mitchell2003two} model for the ionic potential, given by
\begin{equation}
 \begin{aligned}
    &\partial_t v = \frac{1}{\tau_{in}}h v^2 (1-v) - \frac{v}{\tau_{out}}, \\
    &\partial_t h =\left\{ \begin{aligned} 
    &\frac{1-h}{\tau_{open}}\, \text{ is } v<v_{gate},\\
    &\frac{-v}{\tau_{close}}\, \text{ otherwise}.
                    \end{aligned}\right.
\end{aligned}
\label{MS}
\end{equation} Concerning the numerical schemes, the 2D Discrete Duality Finite Volume (DDFV) \cite{delcourte2005discrete} scheme was used for the space discretization and an order 2 Semi-implicit Backward Discretization scheme (SBDF2) for the temporal discretization. Thus, in the bidomain code, the term $\partial_t v + f(v)$ is approached by
\begin{multline}
    \label{sbdf2}
     \frac{1}{dt} \left ( 1.5 v^{n} - 2 v^{n-1} + 0.5 v^{n-2}  \right ) \\
     + 2f(v^{n-1}, h^{n-1})-f(v^{n-2}, h^{n-2}),
\end{multline} for a time step $dt > 0$.
Two 2D meshes were used, an artificial mesh of a circular heart in a circular torso, and a 2D slice of a segmentation of a patient's ventricle and torso.

\begin{figure}[h!]
\centering
\includegraphics[width=6cm]{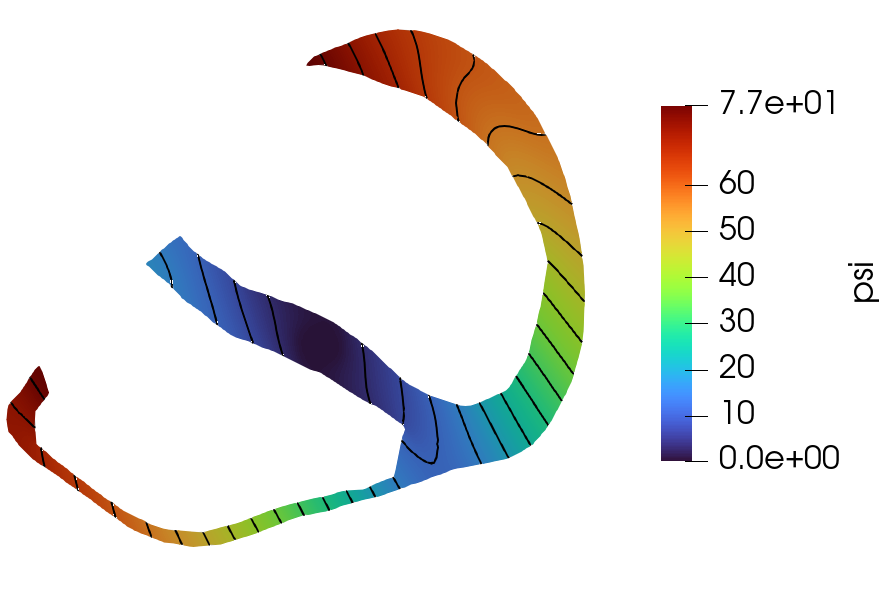}
\caption{Activation times computed from the bidomain propagation.}
\label{fig:ATs}
\end{figure}

The activation map of the heart was computed as a piecewise affine function given by its values $\psi_i$ at each node $x_i$ of the mesh (example Figure \ref{fig:ATs}). The time $\psi_i(x_i)$ is defined as the fist time at which $t \longrightarrow v(x_i)$ crosses a threshold. This threshold is set at $0.5$, as in the Mitchell-Shaeffer model $v$ takes values between $0$ and $0.95$. To find $\psi_i$ at each $x_i$, the functions $t \longrightarrow v(x_i)$ are interpolated by cubic splines.

\subsection{Numerical resolution of the two approached formulations}
The equations \eqref{SF} and \eqref{BF} being static, we used again a 2D DDFV spatial scheme for the resolution and the discretization. Two predefined shapes $\bar{V}$ were designed. The first one is a Heaviside front, smoothed such that the right-hand sides (RHS) of the two formulations $ \partial_t \tilde{v} + f(\tilde{v}) $ and $-\text{div}(\sigma_i \nabla \tilde{v})$ are $L^2$ regular. The second one is the front obtained by solving the 0D Mitchell-Shaeffer problem, with a $\mathcal{C}^{\infty}$ initial stimulation.
For the formulation \eqref{BF}, the RHS is computed by applying a DDFV scheme to the chosen form of $\tilde{v}$. For the source formulation \eqref{SF}, the temporal derivative term is either explicitly given, $\partial_t \tilde{v} = \bar{V}'$, either derived numerically as in \eqref{sbdf2}. In the Source Formulation \eqref{SF}, the ionic term is a cubic function of $v$.

\subsubsection{Choice and computation of the front $\bar{V}$} 
In the first case, $\bar{V}$ is a smoothed Heaviside function depending on a front duration parameter $\varepsilon$, specifically, $\tilde{v} = v_\varepsilon(x, t) = \tilde{H}_\varepsilon(t-\psi(x))$. To fulfill the constraint of a $H^1$ regularity for the solution of \eqref{SF}, we imposed $\tilde{H}_\varepsilon(x) = 0$ for $ x<-\varepsilon$, $\tilde{H}_\varepsilon(x) = 0.94 $ for $ x>\varepsilon$ and $\tilde{H}_\varepsilon(x) = g $ in between ($ |x| \leq \varepsilon$), with $
g(x) = 0.94 \left (\frac{-1}{4 \varepsilon^3}x^3 + \frac{3}{4 \varepsilon} x + \frac{1}{2} \right )$.
Both $v_\varepsilon$ and its temporal derivative are then analytically known.

For the 0D Mitchell-Shaeffer front (MS0D), the Mitchell-Shaeffer model was solved in 0D, using classical python integration tools, on a time range of 330 ms (AP duration). Then, a cubic spline interpolation from python scipy library was used to compute a model for $\tilde{v}$ callable at all times. The time scale was translated to center $\tilde{v}$ on the activation time.

\subsubsection{Computation of the ionic term}
In the Mitchell-Shaeffer model, for an ideal propagation, $h=1$ on the activation front.
Then the Mistchell-Shaeffer fonction $f_\text{MS}$ becomes $f_\text{MS}(v)=\frac{v^2 (1-v)}{\tau_{in}} - \frac{v}{\tau_{out}},$ with $\tau_{in}=0.3$ and $\tau_{out}=6$ as in the original paper \cite{mitchell2003two}. Denoting $f_\text{ref}(v)$ the bidomain ionic term, we compared $f_\text{MS}$ to a cubic polynomial interpolation of $v_{\text{ref}} \longrightarrow f_{\text{ref}}(v_{\text{ref}})$ denoted $f_\text{int}$. Several random points $x_k$ of the heart mesh were selected, and a least-square minimization was realised to find $v(x_k) \longrightarrow f_\text{int}(v(x_k))$ under the constraint that $0$ and $0.94$ are zeros of the function. This means that two coefficients were left to find: the dominant coefficient and the third root. For all points $x_k$, the interpolated coefficients found were close to each other. Thus we took their mean value for designing the function $f_\text{int}(v)$.
\begin{figure}[h!]
\includegraphics[width=7.1cm]{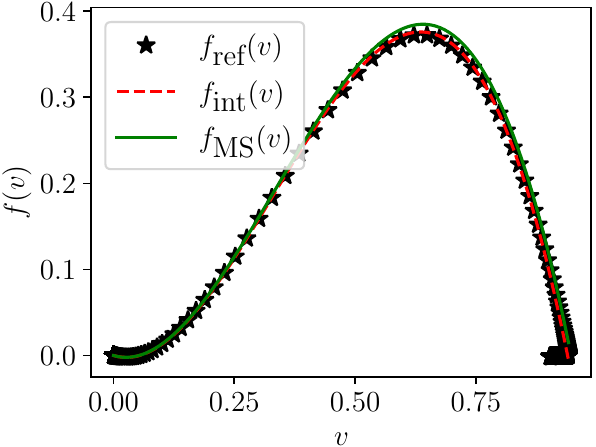}
\caption{Comparison between $f_{\text{MS}}$ in green and $f_{\text{int}}$ in red, as functions of $v$. The stars are the reference values $f_\text{ref}(v)$ obtained by the bidomain propagation.}
\label{fig:interp_f}
\end{figure}As illustrated on one mesh point in Figure \ref{fig:interp_f} the $f_\text{int}$ function fits better the reference curve $f_{\text{ref}}(v_{\text{ref}})$ than $f_\text{MS}$ and was thus chosen for the following.

\section{ Results}

\subsection{Verification of the resolution of the two formulations}
To verify the resolution, we substituted $v_\text{ref}$ in both formulations. For the formulation \eqref{SF}, the RHS is computed with \eqref{sbdf2}. $u_1$ and $u_2$ are expected to match exactly the reference $u_\text{ref}$. In practice, the differences are $\approx 10^{-13}$ due to numerical errors. However, the result $u_1$ is very sensitive to small errors in the RHS. For instance, if we replace formula \eqref{sbdf2} by another discretization scheme of same order, the errors $\frac{\|u_1-u_\text{ref}\|_{L^2(\Omega_T)}}{\|u_\text{ref}\|_{L^2(\Omega_T)}}$ jump to an order $10^{-2}$ or even $10^{-1}$, as presented in Tables \ref{table:1} and \ref{table:2} for the circular mesh.
\begin{table}[htbp]
\label{table:1}
\caption{\label{table:1} Relative $L^2(\Omega_T)$ errors for the solution of \eqref{SF} when changing $f$ in \eqref{sbdf2}.}
\vspace{4 mm}
\centerline{\begin{tabular}{|c|c|} \hline
Ionic term    & $\frac{\|u_1-u_\text{ref}\|_{L^2(\Omega_T)}}{\|u_\text{ref}\|_{L^2(\Omega_T)}}$ \\ \hline
$f_\text{ref} \left (v_\text{ref} \right )$ (reference)   & 1.5E-13 \\
$f_\text{int} \left (v_\text{ref}(t^n) \right )$   & 8.1E-01 \\
$f_\text{MS}\left (h(t^n), v_\text{ref}(t^n) \right )$ & 4.4E-02 \\
\hline
\end{tabular}}
\end{table}
\vspace{-4 mm}
\begin{table}[htbp]
\label{table:2}
\caption{\label{table:2} Relative $L^2(\Omega_T)$ errors for the solution of \eqref{SF} when changing $\partial_t v$ in \eqref{sbdf2}.}
\vspace{4 mm}
\centerline{\begin{tabular}{|c|c|} \hline
Discrete derivative scheme    & $\frac{\|u_1-u_\text{ref}\|_{L^2(\Omega_T)}}{\|u_\text{ref}\|_{L^2(\Omega_T)}}$  \\ \hline
SBDF2 (reference)   & 1.5E-13 \\
Euler centered & 2.5E-02 \\
Explicit Euler & 4.3E-01\\
\hline
\end{tabular}}
\end{table} 
\begin{figure}[h]
\centering
\includegraphics[width=8cm]{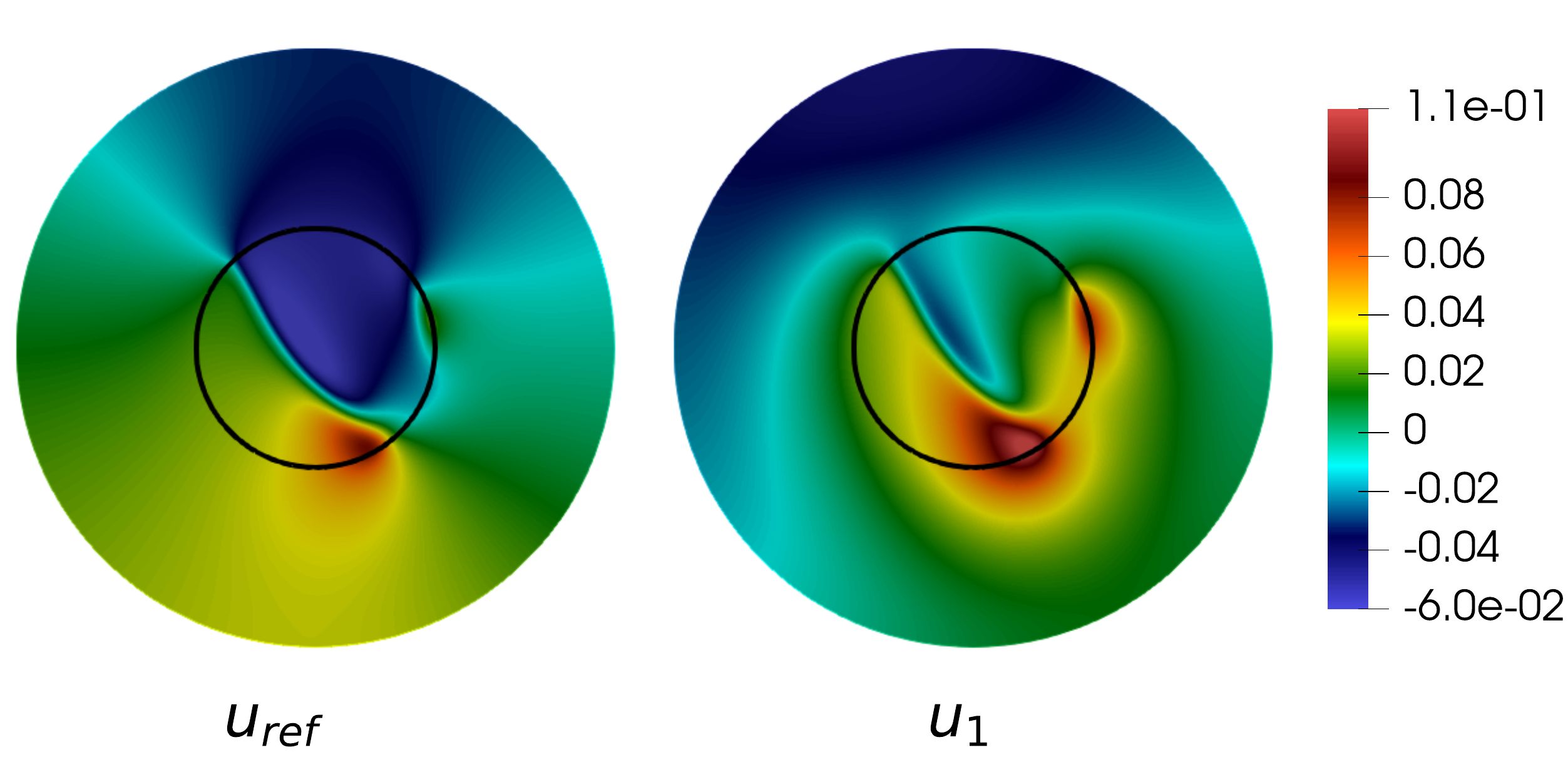}
\caption{Bidomain solution $u_\text{ref}$ on the left, and $u_1$ solution of \eqref{SF} on the right for $\partial_t v = \partial_t v_\text{ref}$ and $f(v) = f_\text{int}\left (v_\text{ref}(t^n) \right )$ in \eqref{sbdf2} (line 2 of Table \ref{table:1}).}
\label{fig:F1 variability ref}
\end{figure}
These results anticipated mediocre performances for \eqref{SF} in case of a pre-shaped front, the source term being very sensitive to small variations, and relying on a balance between the transmembrane voltage time derivative and the ionic term. Figure \ref{fig:F1 variability ref} illustrates this observation.

\subsection{Analysis of the effects of the formulation} 
\begin{figure}[h!]
\centering
\includegraphics[width=8.2cm]{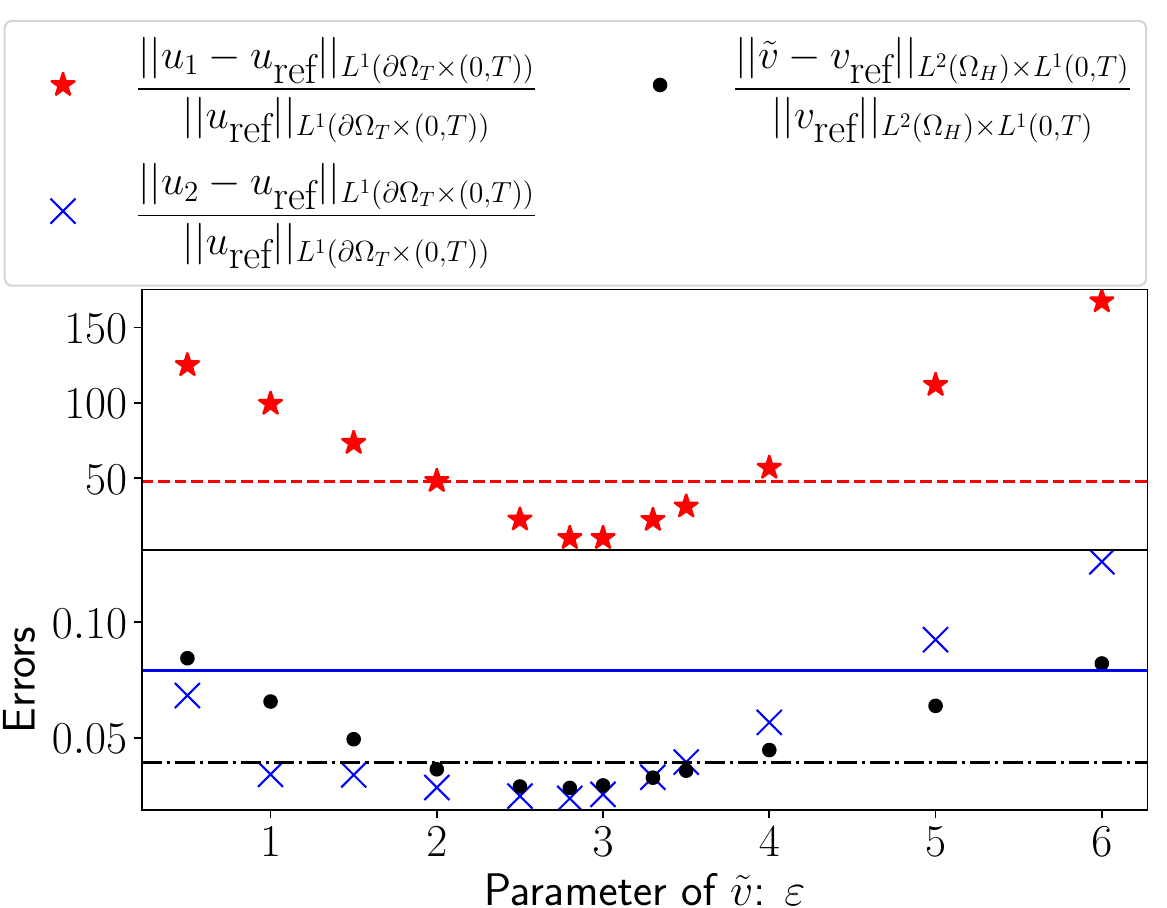}
\caption{Relative differences between $u_1$ (red), $u_2$ (blue) and $u_\text{ref}$ respectively, and $\tilde{v}$ and $v_\text{ref}$ (black). $\tilde{v}$ may be the MS0D model (horizontal dashed lines), or $v_\varepsilon$ with varying duration parameter $\varepsilon$ (points).}
\label{fig:errors epsilon all times}
\end{figure}

In the following, $v$ is substituted by $\tilde{v} = v_\varepsilon$ or the MS0D front. The relative differences between $u_1$ and $u_\text{ref}$ and $u_2$ and $u_\text{ref}$ are presented in Figure \ref{fig:errors epsilon all times}. 
Accordingly to the previous results, \eqref{SF} generates huge errors, thus cannot be used in computations. On the contrary, \eqref{BF} seems robust and generates less than $10\%$ of difference for all the tested $\tilde{v}$. For parameter $\varepsilon$ between $0.5$ and $5$ for $v_\varepsilon$, the errors are lower than the ones obtained with the MS0D model. We can also notice that there exist an optimal $\varepsilon_0$ parameter that minimizes the difference $\frac{||u_2-u_{\textnormal{ref}}||_{L^1\left (\partial \Omega_T \times (0, T)\right )}}{||u_{\textnormal{ref}}||_{L^1\left (\partial \Omega_T \times (0, T)\right )}} $. Anyway, if we look at the difference in space only, $\varepsilon_0$ slightly depends on time.

\subsection{First sensitivity analysis}
Focusing on \eqref{BF}, a white noise perturbation $w$ was added to the RHS of the equation. Four noisy transmembrane voltages were designed: $\tilde{v} = v_\text{ref} + w(x)$, $\tilde{v} = v_\varepsilon + w(x)$, $\tilde{v} = H_\varepsilon(t-\psi(x) + w(x))$ and $\tilde{v} = H_\varepsilon(t-\psi(x) + w_0)$, where $w_0$ is a scalar. The last case represent a displacement of the whole activation front line, as it would be if the initial activation site were uncertain. For the \% of noise on the activation time, a reference amplitude was needed. Denoting $p$ the position of the stimulation site, we allow $p + \Delta p$ to be located in all septum. Then $\Delta t = \frac{\Delta p}{c}$ with $c$ the mean speed of the activation front in the septum.
\begin{figure}[h!]
\centering
\includegraphics[width=8cm]{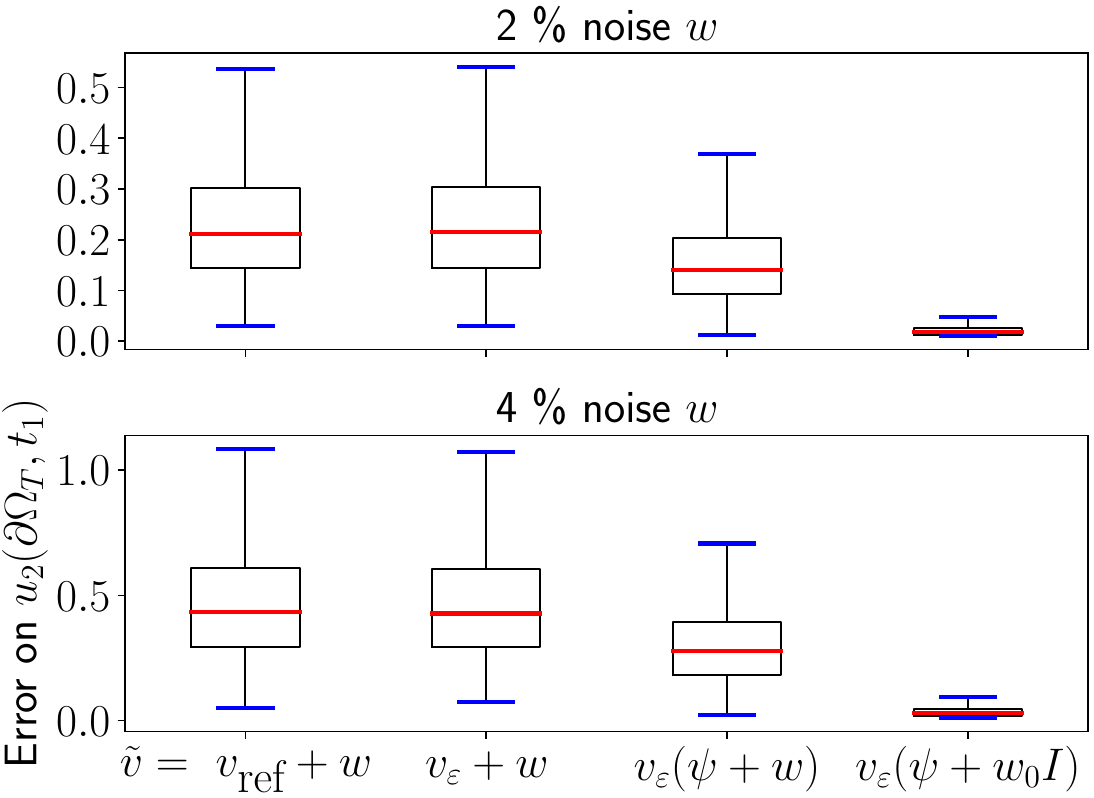}
\caption{Relative difference between $u_2$ and $u_\text{ref}$ on the torso at time $t=35$ ms for $5000$ realisations of white noise $w$. $\varepsilon=2.5$.}
\label{fig:moustache noise t35}
\end{figure}
Surprisingly, the chart Figure \ref{fig:moustache noise t35} shows no significant difference between the solutions of \eqref{BF} with $\tilde{v} = v_\text{ref} + w$ and $\tilde{v} = v_\varepsilon + w$, confirming again the validity of using \eqref{BF} in future computations. Moreover, the formulation does not appear to be very sensitive to small noise in the stimulation site. In this case, the errors are significantly inferior to raw gaussian noise on the activation time. 

\section{Discussion and Conclusion}     
Our results point out clearly that the Source Formulation \eqref{SF} cannot be used for propagating $u$ from $v$. Even with very accurate entries, the body surface potential map is no longer preserved. This emphasizes the underlying equilibrium existing between the ionic term and $v$. On the contrary, the Balance Formulation allows to correctly recover body surface potentials even with a simple preshaped transmembrane potential or with noise. \\
It would be interesting to push further this study in 3D. Moreover, a real sensitivity analysis of $u_2$ as function of $\psi$ would be very useful. In future work, we could also study the sensitivity of formulation \eqref{BF} to heterogeneities in the heart conductivity.

\balance

\section*{Acknowledgments}  
%
This study received financial support from the French Government as part of the “Investments of the Future” program managed by the National Research Agency (ANR), Grant reference ANR-10-IAHU-04.

\bibliography{refs}


  
  
      

\begin{correspondence}
Emma Lagracie\\
200 avenue de la vieille tour, 33400 Talence, France\\
emma.lagracie@inria.fr
\end{correspondence}

\end{document}